\documentclass{amsart}
\usepackage{graphicx}

\numberwithin{equation}{section}

\theoremstyle{plain}
\newtheorem{theorem}{Theorem}
\theoremstyle{remark}
\newtheorem*{remark}{Remark}
\newtheorem*{acknowledgments}{Acknowledgments}

\newcommand\FIG[3]{\begin{figure}
    \includegraphics[#3]{#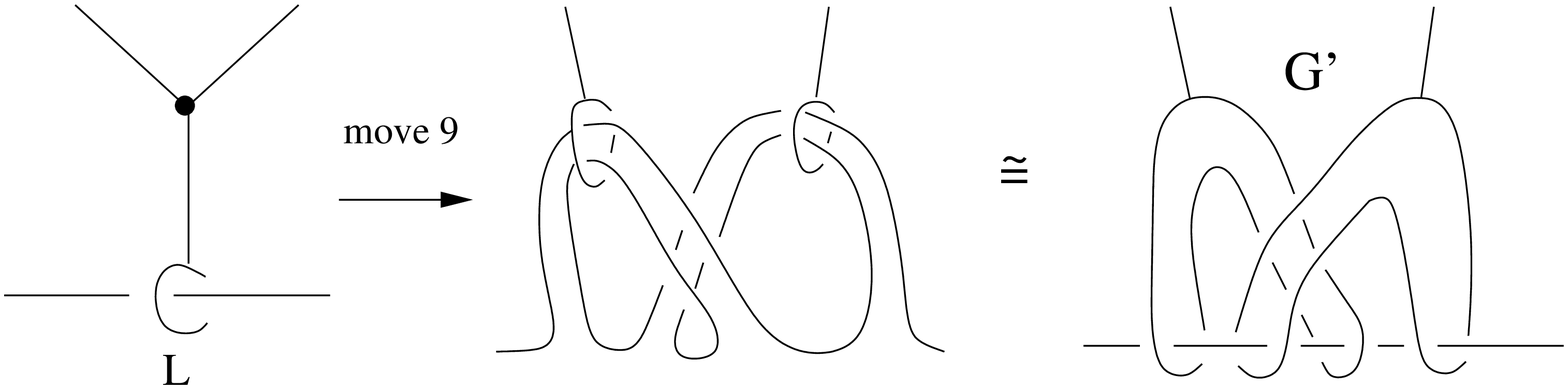}
    \caption{#2}
    \label{fig:#1}
    \end{figure}}

\begin{document}

\title{Replacing a graph clasper by tree claspers}

\author{Kazuo Habiro}

\address{Research Institute for Mathematical Sciences\\Kyoto
  University\\Kyoto\\606-8502\\Japan}

\email{habiro@kurims.kyoto-u.ac.jp}

\begin{abstract}
  We prove that two links related by a surgery along a connected,
  strict graph clasper of degree $n$ are $C_n$-equivalent, i.e,
  related by a sequence of surgeries along strict tree claspers of
  degree $n$.
\end{abstract}

\date{October 21, 2005}

\thanks{This research was partially supported by the Japan Society for
  the Promotion of Science, Grant-in-Aid for Young Scientists (B),
  16740033.}


\maketitle

\section{Introduction}
\label{sec:introduction}

Goussarov \cite{Goussarov:finite,Gusarov:variations} and the author
\cite{H} independently introduced topological calculus of surgery
along claspers.  One of the main achievements in these theories is the
following characterization of the topological information carried by
Goussarov-Vassiliev finite type invariants
\cite{Vassiliev,Gusarov:91-1,Gusarov:91-2,Birman:93,Birman-Lin,Kontsevich,Bar-Natan}:
Two knots in $S^3$ have the same values for any Goussarov-Vassiliev
invariant of degree\;$<n$ if and only if they are related by a
sequence of {\em $C_n$-moves} \cite{Gusarov:variations,H}.  Here a
$C_n$-move is defined as surgery along a certain type of {\em tree
clasper}, which is a framed unitrivalent tree with each univalent
vertex attached to the knot.

In \cite[\S 8.2]{H}, the author also introduced the notion of {\em
graph claspers} for links, which is a generalization of the notion of
tree claspers, where the tree part is replaced by a unitrivalent
graph.  There we explained the idea that graph claspers may be
regarded as {\em topological realizations} of unitrivalent graphs
(also called Feynman diagrams, Jacobi diagrams, etc.) used by
Bar-Natan \cite{Bar-Natan} to describe the structure of the graded
quotients of the Goussarov-Vassiliev filtration.  Recall that, in the
diagram level, any connected, unitrivalent graph diagram on a
$1$-manifold is equivalent under the $STU$ relations to a linear
combination of tree diagrams.  The purpose of this short note is to
prove a topological version of the above-mentioned fact: surgery along
a strict graph clasper $G$ for a link can be replaced by a sequence of
surgeries along strict tree claspers of the same degree as $G$.

\section{Definitions}
\label{sec:definitions}

We freely use the definitions, notations and conventions in \cite{H}.

In the following, $M$ denotes a compact, connected, oriented
$3$-manifold.

A {\em tangle} $\gamma $ in $M$ is a ``link'' in the sense of \cite[\S
1.1]{H}, i.e., a proper embedding $f\colon\thinspace\alpha \rightarrow M$ of a
compact, oriented $1$-manifold $\alpha $ into $M$.  As usual, we
systematically confuse $\gamma $ and the image $\gamma (\alpha )\subset M$.  A link in the
usual sense is a tangle consisting only of circle components.

Two tangles $\gamma $ and $\gamma '$ in $M$ are {\em equivalent}, denoted by
$\gamma \cong\gamma '$, if $\gamma $ and $\gamma '$ are ambient isotopic fixing the
endpoints.

For the definitions of {\em claspers}, {\em tree claspers} and {\em
graph claspers}, see \cite[\S1, \S2, \S8.2]{H}.  Note that a tree
clasper is a special kind of connected graph clasper.  A graph clasper
$G$ is called {\em strict} if $G$ has no leaves.  I.e., a strict graph
clasper is a clasper consisting only of disk-leaves, nodes, and edges.

An important property of a strict graph clasper $G$ for a tangle $\gamma $
in $M$ is that $G$ is {\em tame} (see \cite[\S2.3]{H}), and
consequently surgery along $G$ does not change the $3$-manifold up to
canonical homeomorphism.  (The proof of this fact is similar to
\cite[Proposition 3.3]{H}.)  Thus we may regard the result $\gamma ^G$ from
$\gamma $ of surgery along $G$ as a tangle in $M$.

A disk-leaf $A$ in a clasper for a tangle $\gamma $ is {\em simple} if $A$
intersects $\gamma $ by only one point.  A strict graph clasper is {\em
simple} if all the disk-leaves are simple.

For $k\ge 1$, a {\em $C_n$-move} is a local move on a tangle defined as
surgery along a strict tree clasper of degree $n$.  The {\em
$C_n$-equivalence} on tangles is generated by $C_n$-moves and
equivalence.

\section{Statement and proof of the result}
\label{sec:result}

The purpose of this note is to prove the following.

\begin{theorem}[Stated in a different form in {\cite[\S 8.2, p.68,
	l.4]{H}}]
  \label{r1}
  Let $\gamma $ be a tangle in a compact, connected, oriented $3$-manifold
  $M$, and let $G$ be a strict graph clasper for $\gamma $ in $M$ of degree
  $n\ge 1$, which is not necessarily simple.  Then $\gamma $ and $\gamma ^G$ are
  $C_n$-equivalent.  (Consequently, by \cite[Theorem 3.17]{H}, there
  are finitely many disjoint simple tree claspers $T_1,\ldots ,T_p$ for
  $\gamma $ of degree $n$ such that $\gamma ^G\cong\gamma ^{T_1\cup \dots \cup T_p}$.)
\end{theorem}

\begin{proof}
  We may safely assume that $G$ is connected.

  The proof is by induction on the number $e(G)$ of edges in $G$.  If
  $e(G)=1$, then $G$ is already a strict tree clasper, and hence the
  assertion follows.

  Let $e(G)>1$.  If we have the assertion for the case when $G$ is
  simple, then we have the general case by replacing a single strand
  by a parallel family of strands.  Hence we may assume that $G$ is
  simple.  (This assumption is just for simplifying explanations and
  figures.)  Choose any disk-leaf $L$ of $G$.  Since $e(G)>1$, $L$ is
  joined by an edge to a node.  Let $G'$ denote the the strict graph
  clasper obtained from $G$ by move 9 of \cite[Proposition 2.7]{H} and
  isotopy as depicted in Figure \ref{fig:1}.  We have
  $\gamma ^{G'}\cong\gamma ^G$.  There are two cases.  \FIG{1}{}{height=25mm}

  {\it Case 1.  $G'$ is connected.}\quad Since $e(G')=e(G)-1$ and
  $\deg G'=\deg G$, the assertion follows from the induction
  hypothesis.

  {\it Case 2.  $G'$ consists of two components $G_1$ and $G_2$.}\quad
  We have
  \begin{gather}
    \label{e1}
    \gamma ^{G_1\cup G_2}\cong\gamma ^G,\\
    \label{e2}
    \deg G_1 +\deg G_2 =\deg G =n.
  \end{gather}
  For $i=1,2$, let $N_i$ be a small regular neighborhood of $G_i$,
  such that $N_1\cap N_2$ is empty.  Let $\gamma _i=\gamma \cap N_i$, which is a
  tangle in $N_i$.  Let $c\subset \gamma $ be the component which intersects $L$.
  Let $L_i$ denote the new disk-leaf in $G_i$, which intersects $c$.
  By the induction hypothesis and \cite[Theorem 3.17]{H}, for $i=1,2$,
  there is a clasper $F_i$ consisting of finitely many disjoint,
  simple strict tree claspers of degree $\deg G_i$ for $\gamma _i$ in
  $N_i$, such that
  \begin{equation}
    \label{e3}
    (\gamma _i)^{F_i}\cong(\gamma _i)^{G_i}.
  \end{equation}
  For $i=1,2$, $c\cap N_i$ consists of two components $c_i,c_i'$, where
  these components for $i=1,2$ placed in $c$ in the order
  $c_1,c_2,c_1',c_2'$.  For each of these arcs $c_1,c_2,c_1',c_2'$,
  there are finitely many intersecting disk-leaves and finitely many
  winding edges, as depicted in the left-hand side of Figure
  \ref{fig:2}.  \FIG{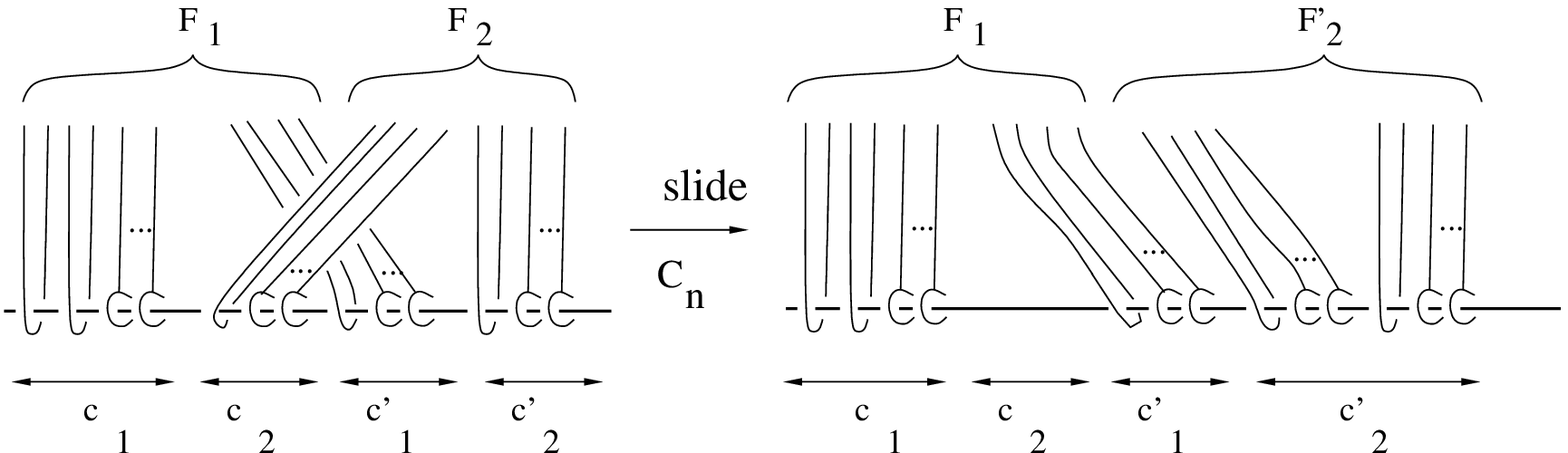}{}{height=35mm} Slide the disk-leaves and edges
  of $F_2$ around $c_2$ along $c$ to traverse those of $F_1$ around
  $c_1'$.  The result is depicted in the right-hand side of Figure
  \ref{fig:2}.  By \eqref{e2} and \cite[Propositions 4.4 and 4.6]{H},
  this sliding does not change the $C_n$-equivalence class of result
  of surgery on $\gamma $.  Let $F_2'\subset N_2$ denote the clasper obtained
  from $F_2$ by the above sliding moves.  It follows from the
  construction of $G_1$ and $G_2$ that $\gamma ^{F_1\cup F_2'}\cong
  \gamma ^{G_1'\cup G_2'}$, where $G_1'$ and $G_2'$ are depicted in Figure
  \ref{fig:3}.  \FIG{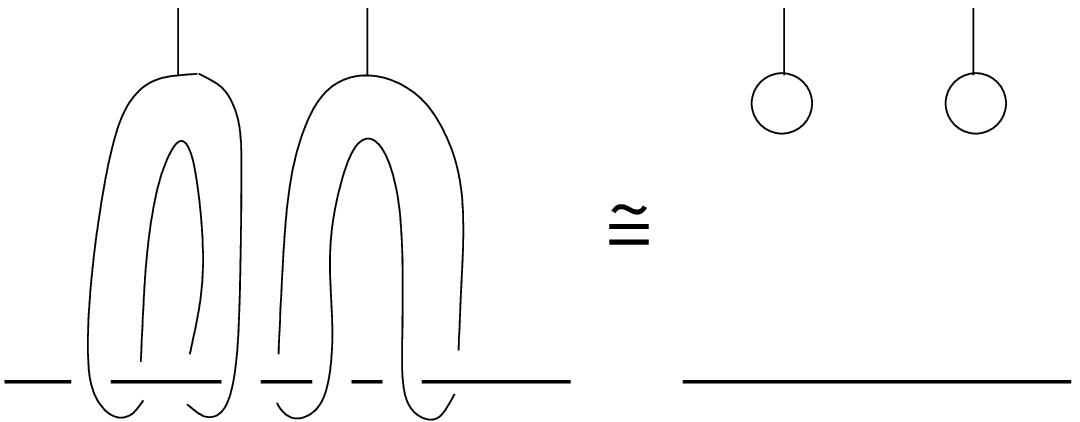}{}{height=25mm} By an obvious graph-clasper
  version of \cite[Proposition 3.4]{H}, we have
  $\gamma ^{G_1'\cup G_2'}\cong\gamma $.  Hence we have
  \begin{equation*}
    \gamma ^{F_1\cup F_2}\underset{C_n}{\sim}\gamma ^{F_1\cup F_2'}
    \cong\gamma ^{G_1'\cup G_2'}
    \cong\gamma .
  \end{equation*}
  The assertion follows from this, \eqref{e1}, and \eqref{e3}.
\end{proof}

\begin{remark}
  More systematic study of graph claspers as announced in \cite[\S
  8.2, \S 8.3]{H} will appear elsewhere.
\end{remark}

\begin{acknowledgments}
  The author thanks Jean-Baptiste Meilhan for reading draft versions
  of this paper and giving me many helpful comments.  He also thank
  Toshifumi Tanaka for asking me about the proof of the theorem, which
  motivated him to write this paper.
\end{acknowledgments}

\end{document}